\documentstyle[11pt]{article} % To be LaTeX'ed with cjp3.sty
\def\ba{\begin{eqnarray}}
\def\ea{\end{eqnarray}}
\def\lb{\label}
\def\ni{\noindent}
\def\nn{\nonumber}

\def\be{\begin{equation}}
\def\ee{\end{equation}}

\def\R{\hat{R}}
\def\id{I}

\def\s{\sigma}

\def\t{\tau}

\def\tr{{\mbox{\rm Tr}\,}}
\def\trq{{\mbox{\rm Tr}_{\mbox{\rm q}\,}}}
\def\Dr{{{\cal D}_r}}
\def\Dl{{{\cal D}_{\ell}}}

\newcommand{\Aut}{\mathop{\mbox{\rm Aut}}\nolimits}

\newcommand{\rank}{\mathop{\mbox{\rm rank}}\nolimits}

\begin{document}
\title{Generalized Cayley-Hamilton-Newton identities}
\author{A. Isaev\\
{\small\em Bogoliubov Laboratory of Theoretical Physics, JINR, 141980  Dubna,}\\
{\small\em
 Moscow region, Russia} \medskip \\
O. Ogievetsky\footnote{On leave of absence from 
P. N. Lebedev Physical Institute, 
Theoretical Department, Leninsky pr. 53, 117924 Moscow, Russia}
~~~and~~~ P. Pyatov\footnote{On leave of absence from
Bogoliubov Laboratory of Theoretical Physics, JINR,
141980  Dubna, Moscow region, Russia}\\
{\small\em Center of Theoretical Physics, Luminy,
13288 Marseille, France}}
\date{}
\maketitle

\begin{abstract} 
\noindent The $q$-generalizations of the two fundamental
statements of matrix algebra -- the Cayley-Hamilton theorem and the 
Newton relations -- to the cases of quantum matrix algebras of
an "RTT-" and of a "Reflection equation" types 
have been obtained in \cite{NT}--\cite{IOPS}.
We construct a family of matrix identities which we call
Cayley-Hamilton-Newton identities and
which underlie the characteristic identity 
as well as the Newton relations 
for the RTT- and Reflection equation algebras, in the sense
that both the characteristic identity and the Newton relations
are direct consequences of the Cayley-Hamilton-Newton identities.
\end{abstract}

\section{Introduction}

Let $V$ be a vector space and $\R\in \Aut (V\otimes V)$
an $\R$-matrix of Hecke type, that is, $\R$
satisfies the Yang-Baxter equation and Hecke condition,
respectively,
\ba
\R_1\R_2\R_1 &=& \R_2\R_1\R_2\;,
\lb{YBE} \\
\quad \R^2 &=& \id + (q-q^{-1}) \R\;.
\label{hecke}
\ea
We use here the matrix notations of \cite{FRT}
(e.g., $\R_1 = \R\otimes\id$, $\R_2 = \id\otimes \R$ in (\ref{YBE}) etc.),
$\id$ is an identity operator and $q\neq 0$ is a numeric parameter.

In this note we deal with quantum matrix
algebras of two types: an RTT-algebra and a 
Reflection equation (RE) algebra. They are
associative unital algebras generated, respectively,
by elements of "$q$-matrices" 
$T = || T^i_j||_{i,j=1,\dots,\dim V}$ and
$L = || L^i_j||_{i,j=1,\dots,\dim V}$
subject to relations
\ba
\R\, T_1T_2 &=& T_1T_2 \, \R\; ,
\lb{rtt}
\\
\lb{rlrl}
\R\, L_1 \R\, L_1 &=& L_1 \R\, L_1 \R \; .
\ea

For both these algebras, 
$q$-versions of the Newton identities and
the Cayley-Hamilton theorem have been recently established
(see \cite{NT}--\cite{IOPS}).
The proofs of these two statements
given  for the $q$-matrix $T$
in \cite{PS2} and \cite{IOPS} 
turn out to be very 
similar ideologically and technically,
which indicates that there should exist
a more wide set of identities containing the
Newton and the characteristic identities
as particular cases. The main object
of the present note is a construction of such
generalized Cayley-Hamilton-Newton (CHN) identities.

We prove a $q$-version of the 
CHN identities for the RTT-algebra case.
The CHN identities for the RE algebra are presented also.
In case when both the RTT- and RE algebras originate
from a quasitriangular Hopf algebra, the CHN identities
for the $q$-matrix $L$ can be derived from those for the $q$-matrix
$T$ by a procedure described in \cite{IOPS}. 
An independent proof of the CHN
identities for the RE algebra will be given elsewhere.

Note that taking $\R =P$, the permutation matrix, one obtains --
from any of the $q$-versions of the CHN theorem --
a set of identities for usual matrices with commuting entries.
It is worth mentioning that the CHN identities
appear to be a new result even for the classical matrix algebra.

\section{Notation}
\setcounter{equation}0

We shall begin with a brief reminder on the $\R$-matrix
technique (a more complete treatment can be found, e.g., in
\cite{G,GPS}).

Assume that $q$ is not a root of unity, that is 
$k_q \equiv (q^k-q^{-k})/(q-q^{-1}) \neq 0$ for any $k=2,3,\dots$~.

Given a Hecke $\R$-matrix, one can construct two series of
projectors, $A^{(k)}$ and $S^{(k)}$,
called $q$-antisymmetrizers and
$q$-symmetrizers, respectively. 
They are defined inductively as
\ba
\lb{antis}
A^{(1)}:=\id\  ,&&
A^{(k)}:={1\over k_q}\,
A^{(k{-}1)}\left(q^{k-1}-(k{-}1)_q\R_{k{-}1}\right)A^{(k{-}1)}\ ,
\\
\lb{simm}
S^{(1)}:=\id\  ,&&
S^{(k)}:= {1\over k_q}\,
S^{(k{-}1)}\left(q^{1-k}+(k{-}1)_q\R_{k{-}1}\right)S^{(k{-}1)}\ .
\ea

Further, assume that the $q$-antisymmetrizers fulfil the conditions
\be
\lb{height}
\rank A^{(n)}=1\ , \quad A^{(n{+}1)}=0
\ee
for some $n$. 
In this case the corresponding $\R$-matrix is called {\em even}
and the number $n$ is called the {\em height} of the $\R$-matrix.

For an $\R$-matrix of finite height $n$ one introduces
the following two matrices 
\ba
\lb{d-a}
\Dr :=  {n_q \over q^n} \tr_{(2 \dots n)} A^{(n)}\ , &&
\qquad \!\!
\Dl :=  {n_q \over q^n} \tr_{(1 \dots n-1)} A^{(n)}\ ,
\ea
Here and below we use notation 
$\tr_{(i_1\dots i_k)}$
to denote the operation of taking traces in the spaces on places
$(i_1\dots i_k)$.

\section{Cayley--Hamilton--Newton identities}
\setcounter{equation}0

Let us consider three sequences of elements 
in the RTT-algebra:
\ba
\lb{s}
s_k(T) \,:= \phantom{q^k} &&\!\!\!\!\!\!\!\!\!\!
\tr_{(1\dots k)}(\R_1\R_2\dots\R_{k-1}T_1 T_2\dots T_k)\ ,
\\
\lb{sig}
\s_k(T) \,:=\, q^k &&\!\!\!\!\!\!\!\!\!\!
\tr_{(1\dots k)}(A^{(k)}T_1 T_2\dots T_k)\ ,
\\
\lb{t}
\t_k(T) \,:=\, q^{-k} &&\!\!\!\!\!\!\!\!\!\! 
\tr_{(1\dots k)}(S^{(k)}T_1 T_2\dots T_k)\ ,
\quad k=1,2,\dots\; .
\ea
We also put $s_0(T) =\s_0(T) = \t_0(T) = 1$.

To clarify the meaning of these elements, consider 
the classical limit $\R = P$. 
Denote $\{x_a\}$ the spectrum of the semisimple part
of an operator $X\in \Aut(V)$.  
Then the elements $s_k(X)$, $\s_k(X)$,
$\t_k(X)$ are  symmetric polynomials in $x_a$. Namely, 
$s_k(X)=\tr X^k =\sum_a x_a^k$ are
{\em power sums}, $\s_k(X)= \sum_{a_1<\dots <a_k}x_{a_1}\dots x_{a_k}$
are {\em elementary symmetric functions}, and
$\t_k(X)= \sum_{a_1\leq\dots\leq a_k}x_{a_1}\dots x_{a_k}$ 
are {\em complete symmetric functions}.
We keep this notation for the elements $s_k(T)$, $\s_k(T)$, $\t_k(T)$ of 
the RTT-algebra also.
  
\vskip 0.2cm 
The $q$-version of power sums $s_k(T)$ has been introduced by
J.M. Maillet, who established their important 
% nonclassical
property
--- the commutativity \cite{M}.
Just as in the classical case, the elementary and complete symmetric functions
admit an expression in terms of the power sums 
(see  Corollary 2 below)
and, hence, the commutativity
property extends to any pair of elements of the sets $\{s_k(T)\}$,
$\{\s_k(T)\}$, $\{\t_k(T)\}$.

\vskip .2cm
If $\R$ is an even $\R$-matrix of height $n$, then  
one has $\s_k(T)=0$ for $k>n$ and
the last nontrivial
element $\s_n(T)$ is proportional
to a quantum determinant of $T$, $det_q T$ 
\be
\lb{qdet}
\s_n(T)=q^n det_q T\ .
\ee
 
Finally, we need an appropriate generalization of the matrix
multiplication in the RTT-algebra.
Inspired by the definition of the quantum power sums (\ref{s}),
one can introduce
two versions of a $k$-th power of the $q$-matrix $T$ \cite{IOPS}:
\ba
\lb{k1}
T^{\underline{k}} &:=& \tr_{(1 \dots k-1)}
\left(\R_1 \R_2 \dots \R_{k-1} T_1 T_2 \dots T_k  \right) \ ,
\\
\lb{k2}
T^{\overline{k}} &:=& \tr_{(2 \dots k)}
\left( \R_1 \R_2 \dots \R_{k-1} T_1 T_2 \dots T_k \right)\ .
\ea
We use the superscripts $\underline{k}$ and $\overline{k}$ here
for denoting different types of the $k$-th power of matrix $T$.
This should not make a confusion with the usual  matrix power
(one has $T^{\underline{k}}=T^{\overline{k}}=T^k$
in the classical limit only).

In the same manner one can introduce a pair of versions of 
{\em $k$-wedge} ({\em $k$-symmetric})
powers of the $q$-matrix $T$, $T^{\underline{\wedge k}}$ and 
$T^{\overline{\wedge k}}$ 
($T^{\underline{{\scriptscriptstyle \cal S}k}}$ and
$T^{\overline{{\scriptscriptstyle \cal S}k}}$),
removing the last or the first trace in the definition
of the elementary (complete) symmetric functions, respectively,
\ba
\lb{wk}
T^{\underline{\wedge k}} := \tr_{(1 \dots k-1)}
\left(A^{(k)} T_1  \dots T_k  \right) &,\quad &
T^{\overline{\wedge k}} := \tr_{(2 \dots k)}
\left( A^{(k)} T_1  \dots T_k \right)\ ,
\\
\lb{sk}
T^{\underline{{\scriptscriptstyle \cal S}k}} := \tr_{(1 \dots k-1)}
\left(S^{(k)} T_1  \dots T_k  \right)  &,\quad &
T^{\overline{{\scriptscriptstyle \cal S}k}} := \tr_{(2 \dots k)}
\left( S^{(k)} T_1  \dots T_k \right)\ .
\ea

\vskip .3cm     
With these definitions we can formulate the main result.
\vskip .2cm     \ni
{\bf Theorem.}
{\bf (Cayley-Hamilton-Newton identities for the RTT-algebra}). \\
{\em
Let $\R$ be Hecke $R$-matrix. For any $j$, the following identities
hold}
\ba
\lb{le}
j_q\, T^{\underline{\wedge j}} &=&
  \sum_{k=0}^{j-1}(-1)^{j-k+1}\s_k(T)\, T^{\underline{j-k}}\ ,
\\
\lb{le2}
j_q\, T^{\overline{\wedge j}} &=&
  \sum_{k=0}^{j-1}(-1)^{j-k+1} T^{\overline{j-k}}\, \s_k(T)\ ,
\\
\lb{le3}
j_q\, T^{\underline{{\scriptscriptstyle \cal S}j}} &=&
  \sum_{k=0}^{j-1}\t_k(T)\, T^{\underline{j-k}}\ ,
\\
\lb{le4}
j_q\, T^{\overline{{\scriptscriptstyle \cal S}j}} &=&
  \sum_{k=0}^{j-1}T^{\overline{j-k}}\, \t_k(T)\ .
\ea
{\bf Proof}. We shall give the details of the proof of
the eq.(\ref{le}).
The relations (\ref{le2})--(\ref{le4}) can be proved analogously.

For $k=1,\dots ,j-1$ we have
\ba
\s_k(T) T^{\underline{j-k}}&=&q^{k}\tr_{(1\dots k)}
(A^{(k)}T_{1}\dots T_k)\tr_{(k+1\dots j-1)}
(R_{k+1}\dots R_{j-1} T_{k+1}\dots T_{j})
\nn
\\
&=&q^k\tr_{(1\dots j-1)}
(A^{(k)}T_1\dots T_k\R_{k+1}\dots\R_{j-1}T_{k+1}\dots T_j)
\nn
\\
&=&q^k\tr_{(1\dots j-1)}
(A^{(k)}\R_{k+1}\dots\R_{j-1} T_{1}\dots T_j)
\lb{vsp}
\ea
We use (\ref{antis}) in the form
$q^kA^{(k)}=(k+1)_qA^{(k+1)}+k_qA^{(k)}\R_k A^{(k)}$
to rewrite (\ref{vsp}) as
\ba
&&(k+1)_q\tr_{(1\dots j-1)}(A^{(k+1)}\R_{k+1}\dots\R_{j-1}T_1\dots T_j)
\nn
\\
&&+ \, k_q\tr_{(1\dots j-1)}(A^{(k)}\R_kA^{(k)}
\R_{k+1}\dots\R_{j-1}T_1\dots T_j)\ .
\nn
\ea
In the last term, the right antisymmetrizer $A^{(k)}$ commutes with
the expression
$R_{k+1}\dots R_{j-1}T_1\dots T_j$,
so one can move $A^{(k)}$ through this expression to the right.
Next, we can move $A^{(k)}$
to the left using the cyclic property of the trace. Finally,
$(A^{(k)})^2=A^{(k)}$ and we obtain
\ba
\s_k(T) T^{\underline{j-k}}&=&
(k+1)_q\tr_{(1\dots j-1)}(A^{(k+1)}\R_{k+1}\dots\R_{j-1}T_1\dots T_j)
\nn
\\
&&+ \, k_q\tr_{(1\dots j-1)}(A^{(k)}\R_{k}\dots\R_{j-1}T_1\dots T_j)\ .
\nn
\ea
We have also $\s_0(T) T^{\underline{j}}\; =T^{\underline{j}}$.
Taking the alternative sum, we obtain the relation (\ref{le}).
\hfill \rule{2.5mm}{2.5mm}

\vskip .3cm \ni
{\bf Corollary 1.} ({\bf Newton identities for the RTT-algebra \cite{PS2}}.) \\
{\em
Let $\R$ be Hecke $R$-matrix.
The following iterative relations hold for the elements 
of the sets $\{s_k(T)\}$, $\{\s_k(T)\}$ and $\{\t_k(T)\}$ 
}
\ba
\label{qNewton}
q^{-j} j_q\,\s_j(T)&=&\sum_{k=1}^{j-1}(-1)^{k-1}
\s_{j-k}(T)\, s_k(T)+(-1)^{j-1}s_j(T)\ ,
\\
\label{qNewton2}
q^j j_q\,\t_j(T)&=&\sum_{k=1}^{j-1}
\t_{j-k}(T)\, s_k(T)+s_j(T)\ ,
\\
\label{qNewton3}
0 &=& \sum_{k=0}^{j}(-1)^k q^{2(j-k)}\t_{j-k}(T)\, \s_k(T)\; ,
\quad \forall j=1,2,\dots\ .
\ea
{\bf Proof.}
To obtain the eqs. (\ref{qNewton}) and (\ref{qNewton2}) one just
takes the last trace (in the space with number $j$) in (\ref{le}) and (\ref{le3}), correspondingly.
The eq. (\ref{qNewton3}) then follows from (\ref{qNewton}) and (\ref{qNewton2}).
\hfill \rule{2.5mm}{2.5mm}

\vskip .3cm \ni
\ni {\bf Corollary 2.}
({\bf Cayley-Hamilton theorem for the RTT-algebra \cite{IOPS}}). \\
{\em
Let $\R$ be even Hecke $\R$-matrix of rank $n$.
The $q$-matrix $T$ satisfies identities
}
\ba
\sum_{k=1}^n \s_{n-k}(T)(-T)^{\underline{k}} + \s_n(T)\, \Dl &=&0\ ,
\lb{hc1}
\\
\sum_{k=1}^n (-T)^{\overline{k}}
\s_{n-k}(T) + \s_n(T)\, \Dr &=&0\ .     
\lb{hc2}
\ea

{\bf Proof.}
Let $j=n$ in (\ref{le}). We also have
$A^{(n)}T_1\dots T_n=A^{(n)}det_q T$. Then, 
the eq. (\ref{hc1}) follows by an application of
(\ref{d-a}) and (\ref{qdet}). 
The eq. (\ref{hc2}) is similarly derived from (\ref{le2}).
\hfill \rule{2.5mm}{2.5mm}

\vskip .3cm \ni
\ni {\bf Corollary 3.}
({\bf Inverse CHN theorem for the RTT-algebra}). \\
{\em
The formulas inverse to the eqs. (\ref{le})--(\ref{le4}) are
\ba
\lb{inv1}
T^{\underline{j}} &=&
\sum_{k=1}^{j}(-1)^{k+1}q^{2(j-k)} k_q\,
\t_{j-k}(T)\, T^{\underline{\wedge k}}\ ,
\\
\lb{inv2}
T^{\overline{j}} &=&
\sum_{k=1}^{j}(-1)^{k+1}q^{2(j-k)}k_q\,
T^{\overline{\wedge k}} \t_{j-k}(T)\ ,
\\
\lb{inv3}
T^{\underline{j}} &=&
\sum_{k=1}^{j} (-1)^{j-k} q^{-2(j-k)} k_q\,
\s_{j-k}(T)\, T^{\underline{{\scriptscriptstyle\cal S}k}}\  ,
\\
\lb{inv4}
T^{\overline{j}} &=&
\sum_{k=1}^{j} (-1)^{j-k} q^{-2(j-k)} k_q\,
T^{\overline{{\scriptscriptstyle\cal S}k}} \s_{j-k}(T)\ .
\ea
}

\ni
{\bf Proof.~}
Consider two lower triangular matrices:
\ba
\nn
H&:=&\{ H^j_k = q^{2(j-k)}\t_{j-k}(T) \mbox{~if~} j\geq k\, ;\;\;\;
H^j_k = 0\ \mbox{ otherwise}\}\, , 
\\
\nn
E&:=&\{E^j_k = (-1)^{j-k}\s_{j-k}(T)  \mbox{~if~} j\geq k\, ;\;
E^j_k = 0\ \mbox{ otherwise}\}\, .
\ea
By the eq. (\ref{qNewton3}) 
one has
$HE = \id$.

With this notation 
one rewrites (\ref{le}) as
$
(-1)^{j+1}j_q\, T^{\underline{\wedge j}} = 
\sum_{k=1}^j E^j_k T^{\underline{k}}\ . 
$
Then
$
T^{\underline{j}} = \sum_{k=1}^j 
(-)^{k+1} k_q H^j_k T^{\underline{\wedge k}}\ ,
$
which is equivalent to (\ref{inv1}).

The relations (\ref{inv2})--(\ref{inv4}) are proved  similarly.
\hfill \rule{2.5mm}{2.5mm}
\vskip .5cm

We conclude by formulating the CHN theorem for the RE
algebra.

\newpage
%\vskip .4cm     
\ni
{\bf Theorem.}
{\bf (Cayley-Hamilton-Newton identities for the RE algebra}). \\
{\em
Let $\R$ be Hecke $R$-matrix and  the $q$-matrix $L$ generate the
RE algebra (\ref{rlrl}).
Then the following identities
hold
\be
\lb{le-re}
j_q\, L^{\wedge j} =
  \sum_{k=0}^{j-1}(-1)^{j-k+1}\s_k(L)\, L^{j-k}\ ,
\qquad
j_q\, L^{{\scriptscriptstyle\cal S} j} =
  \sum_{k=0}^{j-1}\t_k(L)\, L^{j-k}\ .
\ee

Here the notation is as follows:
$$
L^{\wedge k} := \trq_{(2\dots k)}
(A^{(k)}L_{\overline{1}}\dots L_{\overline{k}})\ , \qquad
L^{{\scriptscriptstyle\cal S} k} := \trq_{(2\dots k)}
(S^{(k)}L_{\overline{1}}\dots L_{\overline{k}})
$$
are the $k$-wedge and the $k$-symmetric powers 
of the $q$-matrix $L$, respectively;
$L^k$ is the usual matrix power; $\s_k(L) := q^k \trq L^{\wedge k}$ and
$\t_k(L) := q^{-k} \trq L^{{\scriptscriptstyle\cal S}k}$ are the elementary
and complete symmetric functions on the 
spectrum of $L$, respectively; 
$\trq X := \tr (\Dr X)$ is a $q$-trace operation, and $L_{\overline{k}}$
is defined inductively by
$$
L_{\overline{1}} := L_1\ , \qquad
L_{\overline{k}} := \R_{k-1}L_{\overline{k-1}}\R^{-1}_{k-1}\ .
$$
}
\vskip .2cm \ni 
{\bf Acknowledgements:~} We are grateful to  D. Gurevich and
P. Saponov for discussions.
This work is supported in part by the grants for promotion of 
french--russian scientific cooperation: the CNRS grant PICS No. 608 and
the RFBR grant No. 98-01-2033.  
The work of P.P. and A.I. is also partly supported by
the RFBR grant No. 97-01-01041. 
% and by INTAS grant 93-127-ext. 

\end{document}